\documentclass[12pt]{amsart}

\usepackage{amssymb}

\usepackage{enumerate}

\usepackage{graphicx}

\makeatletter
\@namedef{subjclassname@2020}{%
  \textup{2020} Mathematics Subject Classification}
\makeatother

\newtheorem{thm}{Theorem}[section]

\newtheorem{lem}[thm]{Lemma}

\theoremstyle{definition}

\numberwithin{equation}{section}

\begin{document}


\baselineskip=17pt



\title[]{ SUPERCYCLIC VECTORS OF OPERATORS ON\\ NORMED LINEAR SPACES}

\author[Mohammad Ansari]{Mohammad Ansari}
\date{}
\begin{abstract} We give an affirmative answer to a question asked by Faghih-Ahmadi and Hedayatian regarding supercyclic vectors. We show that
if $\mathcal X$ is an infinite-dimensional normed linear space and $T$ is a supercyclic operator on $\mathcal X$, then for any supercyclic vector $x$ for $T$, there exists 
a strictly increasing sequence $(n_k)_k$ of positive integers such that the closed linear span of the set $\{T^{n_k}x: k\ge 1\}$ is not the whole $\mathcal X$.
\end{abstract}

\maketitle
\section{\bf Introduction}
Let $\mathcal X$ be an infinite-dimensional normed linear space and $B(\mathcal X)$ be the space of all bounded linear operators on $\mathcal X$. 
An operator $T\in B(\mathcal X)$ is said to be {\it supercyclic} if there is some $x\in \mathcal X$ such that the set $\{cT^nx: c\in \Bbb C, n\ge 0\}$ 
is dense in $\mathcal X$. In this case, $x$ is called a {\it supercyclic vector} for $T$. To see detailed information about supercyclic operators we 
refer the readers to \cite{bm} and \cite{gp}.
\par For a subset $M$ of $\mathcal X$, by $\bigvee M$ we mean the closed linear span of the set $M$. It is clear that if $x$ is a supercyclic vector 
for an operator $T\in B(\mathcal X)$, then $\bigvee \{T^nx: n\ge 0\}=\mathcal X$. A natural question which may be asked here is that whether 
there is a strictly increasing sequence $(n_k)_k$ of positive integers such that $\bigvee \{T^{n_k}x: k\ge 1\}\neq \mathcal X$.
 In their recently published paper, Faghih-Ahmadi and Hedayatian have proved the following interesting result.
\begin{thm}[Theorem $1$ of \cite{fh}] Let $\mathcal H$ be an infinite-dimensional Hilbert space. If $x$ is a supercyclic vector for $T\in B(\mathcal H)$, then
there is a (strictly increasing) sequence $(n_k)_k$ of positive integers such that $\bigvee \{T^{n_k}x : k\ge 1\}\neq \mathcal H$.
\end{thm}
Then they have asked whether the assertion is true for locally convex spaces or at least for Banach spaces \cite[Question $1$]{fh}. 
In this note, we answer their question affirmatively for normed linear spaces.
 To prove our result, we use Lemma $1.3$ which is an analogue of the following lemma.
\begin{lem}[Lemma $2.3$ of \cite{bkk}] Let $\mathcal A$ be a dense subset of a Banach space $\mathcal X$ and $e$ be a fixed element with $\|e\|>1$. 
Then, for every finite-dimensional subspace $Y\subset \mathcal X$ with $\text{dist}(e, Y )>1$, for every $\epsilon >0$ and $y\in Y$, there is an $a\in \mathcal A$ such that
$\|y-a\|<\epsilon$ and $\text{dist}(e, \bigvee\{Y, a\})>1$.
\end{lem}
 The proof of the above lemma shows that it can also be stated for infinite-dimensional normed linear spaces. Thus, we can give the following modified version.
\begin{lem} Let $\mathcal A$ be a dense subset of an infinite-dimensional normed linear space $\mathcal X$ and $e\in \mathcal X$ be a fixed element with $\|e\|>1$. 
Then, for every finite-dimensional subspace $Y$ of $\mathcal X$ with $\text{dist}(e,Y)>1$, there is an $a\in \mathcal A$ such that
$\text{dist}(e, \bigvee\{Y, a\})>1$.
\end{lem}
\section{\bf Main result}
Now, we are ready to answer Question $1$ of \cite{fh} for normed linear spaces.
\begin{thm} Let $\mathcal X$ be an infinite-dimensional normed linear space. If $x$ is a supercyclic vector for $T\in B(\mathcal X)$, then
there is a strictly increasing sequence $(n_k)_k$ of positive integers such that $\bigvee \{T^{n_k}x : k\ge 1\}\neq \mathcal X$.
\begin{proof} It is clear that every nonzero scalar multiple of $x$ is also a supercyclic vector for $T$.
On the other hand, since $x$ and $Tx$ are linearly independent vectors, we have $\text{dist}(x,\bigvee\{Tx\})>0$.
Thus, without loss of generality, we can assume that $\|x\|>1$ and $\text{dist}(x,\bigvee\{Tx\})>1$. 
Therefore, in view of Lemma $1.3$, if we put $$\mathcal A=\{cT^nx: c\in \Bbb C, n>1\}, Y=\bigvee\{Tx\},$$ and $e=x$, then there is
some $a=cT^{n_2}x\in \mathcal A$ such that $$\text{dist}(x,\bigvee\{Tx,T^{n_2}x\})=\text{dist}(x, \bigvee\{Y,a\})>1.$$ 
Now, let $$\mathcal A_2=\{cT^nx: c\in \Bbb C, n> n_2\}, Y_2=\bigvee\{Tx,T^{n_2}x\}.$$ Then
there is some $a_2=c_2T^{n_3}x\in \mathcal A_2$ such that $$\text{dist}(x,\bigvee\{Tx,T^{n_2}x,T^{n_3}x\})=\text{dist}(x,\bigvee\{Y_2,a_2\})>1.$$ 
By continuing this construction, assume that for some $k\ge 2$, the dense set $\mathcal A_k$ and the finite-dimensional subspace $Y_k$ have been presented
 and (by using Lemma $1.3$) we have found an element $a_k=c_kT^{n_{k+1}}x\in \mathcal A_k$ such that
 $$\text{dist}(x,\bigvee\{Tx,T^{n_2}x,\cdots, T^{n_{k+1}}x\})=\text{dist}(x,\bigvee\{Y_k,a_k\})>1.$$ 
Then we put $$\mathcal A_{k+1}=\{cT^nx: c\in \Bbb C, n>n_{k+1}\}, Y_{k+1}=\bigvee \{Tx, T^{n_2}x, \cdots, T^{n_{k+1}}x\}.$$
Again, by Lemma $1.3$, there is some $a_{k+1}=c_{k+1}T^{n_{k+2}}x\in \mathcal A_{k+1}$ such that 
$$\text{dist}(x,\bigvee\{Tx,T^{n_2}x,\cdots, T^{n_{k+2}}x\})=\text{dist}(x,\bigvee\{Y_{k+1},a_{k+1}\})>1.$$ 
This inductive procedure gives a strictly increasing sequence $(n_k)_k$ (with $n_1=1$) and it is easily seen that $$\text{dist}(x,\bigvee\{T^{n_k}x: k\ge 1\})\ge 1$$
(see the last paragraph of the proof of Theorem $1$ of \cite{fh}).
This shows that $x\notin \bigvee\{T^{n_k}x: k\ge 1\}$ and we are done.
\end{proof}
\end{thm}
The interested readers are invited to investigate Question $1$ of \cite{fh} for operators on locally convex spaces.

\vspace*{.2cm}
Department of Basic Sciences\\
Azad University of Gachsaran\\
Gachsaran, Iran\\
Email Address: \email{mansari@iaug.ac.ir}, \email{ansari.moh@gmail.com}

\end{document}